\documentclass[12pt,reqno]{amsart}

\usepackage[arrow,matrix,curve]{xy}

\usepackage[dvips]{graphicx} 

\usepackage{amssymb, latexsym, amsmath, amscd, array,
}

\theoremstyle{definition}

\numberwithin{equation}{section}

\newcommand\N {{\mathbb N}} 

\newcommand\R {{\mathbb R}}
\newcommand\Q {{\mathbb Q}}

\newcommand\RRR{\mbox{I\!I\!R}}

\begin{document}

\thispagestyle{empty}

\title{Cauchy's continuum}

\author{Karin U. Katz and Mikhail G. Katz$^{0}$}

\address{Department of Mathematics, Bar Ilan University, Ramat Gan
52900 Israel} \email{katzmik@macs.biu.ac.il}

\footnotetext{Supported by the Israel Science Foundation grant
1294/06}



\begin{abstract}
Cauchy's sum theorem of 1821 has been the subject of rival
interpretations ever since Robinson proposed a novel reading in the
1960s.  Some claim that Cauchy modified the hypothesis of his theorem
in 1853 by introducing uniform convergence, whose traditional
formulation requires a pair of independent variables.  Meanwhile,
Cauchy's hypothesis is formulated in terms of a single variable~$x$,
rather than a pair of variables, and requires the error
term~$r_n=r_n(x)$ to go to zero at {\em all\/} values of~$x$,
including the infinitesimal value generated by~$\frac{1}{n}$,
explicitly specified by Cauchy.  If one wishes to understand Cauchy's
modification/clarification of the hypothesis of the sum theorem in
1853, one has to jettison the automatic translation-to-limits.
\end{abstract}

\maketitle

\tableofcontents

\section{Sifting the chaff from the grain in Lagrange}
\label{intro1}

One of the most influential scientific treatises in Cauchy's era was
J.-L. Lagrange's {\em M\'ecanique Analytique\/}, the second edition of
which came out in 1811, when Cauchy was barely out of his teens.
Lagrange opens his treatise with an unequivocal endorsement of
infinitesimals.  Referring to the system of infinitesimal calculus,
Lagrange writes:
\begin{quote}
Lorsqu'on a bien con\c cu l'esprit de ce syst\`eme, et qu'on s'est
convaincu de l'exactitude de ses r\'esultats par la m\'ethode
g\'eom\'etrique des premi\`eres et derni\`eres raisons, ou par la
m\'ethode analytique des fonctions d\'eriv\'ees, on peut employer les
infiniment petits comme un instrument s\^ur et commode pour abr\'eger
et simplifier les d\'emonstrations \cite[p.~iv]{Lag}.%
\footnote{``Once one has duly captured the spirit of this system
[i.e., infinitesimal calculus], and has convinced oneself of the
correctness of its results by means of the geometric method of the
prime and ultimate ratios, or by means of the analytic method of
derivatives, one can then exploit the infinitely small as a reliable
and convenient tool so as to shorten and simplify proofs''.}
\end{quote}
Lagrange's renewed enthusiasm for infinitesimals in 1811 went
hand-in-hand with a reliance both on his method of power series, and
on the principle of the ``generality of algebra'' which proved to be
extremely fruitful throughout the 18th century.  However, Cauchy was
already becoming aware of the limitations of these techniques.  He was
aware of examples such as~$e^{-1/x^2}$ where the Taylor series at the
origin does not reproduce the function; the use of divergent power
series was recognized as leading to errors; the limitations of the
``generality of algebra'' were beginning to be felt, particularly in
the study of Fourier series.  The tension resided in the challenge
posed by Lagrange's treatise: can Cauchy sift the chaff from the
grain?  By~1823, Cauchy was ready to go on the offensive, explicitly
naming the {\em M\'ecanique analytique\/} as the target of his
criticisms.  Cauchy's great accomplishment was his recognition that,
while Lagrange's flawed power series method and his principle of the
generality of algebra do not measure up to the standard of rigor
Cauchy sought to uphold in his own work, the infinitesimals can indeed
be {\em reconciled\/} with such a standard of rigor. The resolution of
the tension between the rejection of Lagrange's conceptual framework,
on the one hand, and the acceptance of his infinitesimals, on the
other, is expressed by Cauchy in the following terms:
\begin{quote}
My main aim has been to {\em reconcile\/} the rigor, which I have made
a law in my {\em Cours d'Analyse\/}, with the simplicity that comes
from the direct consideration of infinitely small quantities (Cauchy
1823, see \cite[p.~10]{Ca23}) [emphasis added--authors].
\end{quote}
Cauchy reconciled his rejection of Lagrange's flawed conceptual
framework, on the one hand, with his acceptance of Lagrange's
infinitesimals, on the other.%
\footnote{We therefore reject Schubring's concept of a Cauchyan
``compromise'' whereby infinitesimals allegedly ``disagree''
\cite[p.~439]{Sch} with rigor.}
The {\em Cours d'Analyse\/} is Cauchy's 1821 textbook~\cite{Ca21}
where his infinitesimal definition of continuity first appeared, as
discussed in the next section.

\section{Cauchy's continuity}
\label{two}

In 1853, at the dusk of his scientific carrier, Cauchy reaffirmed the
definition of continuity he had given 32 years earlier, in his
influential textbook {\em Cours d'Analyse\/}, in the following terms:
\begin{quote}
In accordance with the definition proposed in my {\em Analyse
Alg\'ebrique\/}, and generally accepted today, a function~$u$ of a
real variable~$x$ will be {\em continuous\/} between two given bounds
of~$x$, if [...] an {\em infinitely small\/} increment of the variable
always produces, between the given bounds, an {\em infinitely small\/}
increment of the function itself \cite[p.~32]{Ca53} [emphasis
added---authors].
\end{quote}
Meanwhile, infinitesimals themselves are defined in terms of {\em
variable quantities\/} becoming arbitrarily small (which have often
been interpreted as null sequences).  Cauchy writes that such a null
sequence ``becomes'' an infinitesimal~$\alpha$.  Cauchy's terminology
was analyzed by Sad {\em et al\/} \cite{STB}.  It is interesting to
note that Cauchy suppresses the index of the~$n$-th term in such a
sequence, and emphasizes the competing index of the {\em
order\/},~$n$, of the infinitesimal~$\alpha^n$; this is dealt with in
more detail in Borovik \& Katz \cite{BK}.

We will return to Cauchy's 1853 article below.  In a recent article
attempting a synthesis of Lakoff and Lakatos, T.~Koetsier writes:
\begin{quote}
In the following reconstruction I will interpret some of Cauchy's
results in accordance with the traditional view of his work.  A~[...]
presentation of this view is in [J.]~Grabiner (1981)
\cite[footnote~13]{Ko}.
\end{quote}

What is the {\em traditional view\/}, in accordance with which
Koetsier seeks to {\em interpret\/} some of Cauchy's results?
Briefly, while Cauchy's definition of continuity is ostensibly
formulated in the language of infinitely small quantities, the
traditional interpretation seeks to subsume Cauchy's infinitesimals
under a notion of limit.%
\footnote{The traditional approach to Cauchy was critically analyzed
in a 1973 text by Hourya Benis Sinaceur~\cite{Si}.}

The proclivity to interpreting infinitesimals as limits is hardly
limited to Cauchy's work.  J.~Havenel \cite{Ha} describes the reaction
to such an interpretation, on the part of the American philosopher
C.~S.~Peirce.  Havenel notes that
\begin{quote}
Peirce was fully aware that in his time, the notion of infinitesimal
was strongly rejected by most mathematicians, especially in analysis,
with the works of Weierstrass, Dedekind, and Cantor \cite[p.~101]{Ha}.
\end{quote}
Peirce wrote that
\begin{quote}
the doctrine of limits has been invented to evade the difficulty, or
according to some as an exposition of the significance of[,] the word
{\em infinitesimal\/} \cite[3.122]{NEM} [emphasis added---authors].
\end{quote}

Thus, the traditional interpretation institutes a kind of an automated
infinitesimal-to-limits translation.  Such an interpretation, as
applied to Cauchy's work, is actually considerably older than
Grabiner's book referred to by Koetsier.  Boyer \cite{Boy} had already
declared that Cauchy's infinitesimals ``are to be understood in terms
of limits''.%
\footnote{See our Section~\ref{boyer} below for an analysis of Boyer's
views.}
Meanwhile, the Cauchy scholar P.~Jourdain in his seminal 1913
work~\cite{Jo} makes no attempt to reduce Cauchy's infinitesimals to
limits.

Nearly a century after Jourdain, the limit interpretation has become
so entrenched as to be taken as the literal meaning of Cauchy's
definitions by a number of historians.  Thus, J.~Gray lists {\em
continuity\/} among concepts Cauchy allegedly defined
\begin{quote}
using careful, if not altogether unambiguous, {\em limiting\/}
arguments \cite[p.~62]{Gray08} [emphasis added--authors].
\end{quote}
Similarly, in his 2007 anthology \cite{Haw}, S.~Hawking reproduces
Cauchy's {\em infinitesimal\/} definition of continuity on page
639--but claims {\em on the same page\/}, in a comic {\em
non-sequitur\/}, that Cauchy ``was particularly concerned to banish
infinitesimals''.

The subject of Cauchy's continuity (and particularly the related sum
theorem) was recently taken up by K.~Br\aa ting \cite{Br}.  We will
build on her work, not so much to restore Cauchy's infinitesimals to
their rightful place in Cauchy's infinitesimal-enriched continuum, as
to argue that the traditional interpretation in terms of limits in the
context of a standard Archimedean continuum, is self-contradictory,
and in particular untenable.

\section{Br\aa ting's close reading}
\label{intro}

The precise relation between Cauchy's variable quantities, on the one
hand, and his infinitesimals, on the other, has been the subject of an
ongoing debate.  In K.~Br\aa ting's text~\cite{Br}, the two sides of
the debate are represented by E.~Giusti \cite{Giu} and D.~Laugwitz
\cite{Lau87}.  Their respective positions can be traced to a pair of
rival interpretations of the continuum found in the work of
E.~Bj\"orling, a contemporary of Cauchy.

In a footnote to his 1846 paper \cite{Bj46}, Bj\"orling for the first
time introduces his distinction between the following two clauses:
\begin{enumerate}
\item[(A)]
\label{bj2}
``for every {\em given\/} value of~$x$'';
\item[(B)]
\label{bj1}
``for all values of~$x$''.
\end{enumerate}
Here clause~(A) refers to what we would describe today as the ``thin''
real Archimedean continuum, or A-continuum.  Meanwhile, clause~(B)
describes the broader class, including elements described by
Bj\"orling as variable quantities, more concretely sequences depending
on~$n$, corresponding to an enriched Leibnizian continuum.  Johann
Bernoulli was the first to use infinitesimals systematically as a
foundational concept.  Therefore we will refer to such a ``thick''
continuum as a Bernoullian continuum, or B-continuum.%
\footnote{\label{lawvere}An infinitesimal-enriched B-continuum is not
a unique mathematical structure.  Thus, the intuitionistic
Nieuwentijdt-Lawvere continuum is a markedly different implementation
of an infinitesimal-enriched continuum as compared to Robinson's, as
it contains nilsquare infinitesimals, see J. Bell \cite{Bel} and
Appendix~\ref{rival2} below.  See Feferman~\cite{Fef} for an analysis
of the continuum in terms of predicativism and conceptual
structuralism.}
A more detailed discussion of the rival continua may be found in
Appendices~\ref{rival1} and~\ref{rival2}.

Terminology similar to Bj\"orling's was exploited by S.~D.~Poisson.
Poisson describes infinitesimals as being ``less than any {\em
given\/} magnitude of the same nature'' [emphasis added--authors].%
\footnote{Quote from Poisson~\cite[p.~13-14]{Po} reproduced in Boyer
\cite[p.~283]{Boy}.  Note that P.~Ehrlich inexplicably omits the
crucial modifier ``given'' when quoting Poisson in footnote~133 on
page~76 of \cite{Eh}.  Based on the incomplete quote, Ehrlich proceeds
to agree with Veronese's assessment (of Poisson) that "[t]his
proposition evidently contains a contradiction in terms"
\cite[p.~622]{Ve91}.  Our assessment is that Poisson's definition of
infinitesimals is consistent if understood in terms of Bj\"orling's
dichotomy.}
The distinction between {\em given\/}, i.e.~constant, as opposed to
{\em variable\/}, i.e.~changing, is close to Bj\"orling's dichotomy.

Br\aa ting was hardly the first to analyze the fundamental difference
between the two continua.  Having outlined the developments in real
analysis associated with Weierstrass and his followers, Felix Klein
pointed out in 1908 that
\begin{quote}
The scientific mathematics of today is built upon the series of
developments which we have been outlining.  But {\em an essentially
different conception of infinitesimal calculus has been running
parallel with this [conception] through the centuries\/}
\cite[p.~214]{Kl} [emphasis added---authors].
\end{quote}
Thus we have two parallel tracks for conceptualizing infinitesimal
calculus:
\[
\xymatrix@C=95pt{{} \ar@{-}[rr] \ar@{-}@<-0.5pt>[rr]
\ar@{-}@<0.5pt>[rr] & {} 
& \hbox{\quad
B-continuum} \\ {} \ar@{-}[rr] & {} & \hbox{\quad A-continuum} }
\]
(this theme is pursued further in Appendix~\ref{rival1}).  Klein
further points out that such a parallel conception of calculus
\begin{quote}
harks back to old metaphysical speculations concerning the {\em
structure of the continuum\/} according to which this was made up of
[...] infinitely small parts \cite[p.~214]{Kl} [emphasis
added---authors].
\end{quote}
The rival theories of the continuum evoked by Klein are the subject of
Bj\"orling's deliberations here, as well.

In his 1853 text \cite{Bj53}, Bj\"orling exploits this distinction to
argue against a purported counterexample, published by F.~Arndt
\cite{Ar} in 1852, to Cauchy's 1821 ``sum theorem''.%
\footnote{Br\aa ting \cite[p.~521]{Br} translates Cauchy's sum theorem
as follows: ``When the different terms of the series~$[u_0 + u_1 + u_2
+ \cdots + u_n + \cdots ]$ are functions of the same variable~$x$,
continuous with respect to that variable in the vicinity of a
particular value for which the series is convergent, the sum~$s$ of
the series is also a continuous function of~$x$ in the vicinity of
this particular value.''  The reference in Cauchy
is~\cite[p.~131-132]{Ca21}.}
Namely, Bj\"orling points out that in fact Arndt's counterexample only
converges ``for every {\em given\/} value'', i.e., value from the
narrow A-continuum.  Meanwhile, it does not converge ``for all
values'', i.e., values from the enriched B-continuum.  Bj\"orling
concludes that Cauchy's 1853 hypothesis in fact bars Arndt's example.%
\footnote{Br\aa ting does not comment on how the hypothesis of
Cauchy's 1821 sum theorem may have been viewed by Bj\"orling.}

The mutual interactions and influences between Cauchy and Bj\"orling
were explored by Grattan-Guinness \cite{Grat86}, who argues that
Cauchy read Bj\"orling's text, and was influenced by it to
modify/clarify the hypothesis of the 1821 sum theorem.  Namely, in
1853 Cauchy added the word {\em always\/} to indicate that the
hypothesis is interpreted to apply for {\em all\/}~$x$ (B-continuum)
rather then merely for every {\em given\/}~$x$ (A-continuum).

Whether Cauchy's addition, in 1853, of the word ``always'' is a {\em
modification\/} or a {\em clarification\/} of the 1821 condition, is
subject to dispute, and is not a major concern here.  A narrow
A-continuum interpretation of the 1821 hypothesis (which would then
falsify the ``sum theorem'' as stated in 1821) is consistent with
Grattan-Guinness's view that Cauchy was influenced by Bj\"orling in
1853 to {\em broaden\/} the interpretation to a~B-continuum.
Laugwitz~\cite[p.~265]{Lau87} quotes Cauchy \cite[p.~31-32]{Ca53} as
admitting that the {\em statement\/} of the 1821 theorem (but not its
proof) was incorrect: ``Au reste, il est facile de voir comment on
doit modifier {\em l'\'enonc\'e\/} du th\'eor\`eme, pour qu'il n'y
plus lieu \`a aucune exception''.  Note that only a single independent
variable,~$x$, occurs in Cauchy's hypothesis, whether in 1821 or in
1853.  As traditionally stated, uniform convergence is a {\em
global\/} condition stated in terms of a {\em pair\/} of independent
variables.  Interpreting Cauchy's addition of the word {\em always\/}
as ``strengthening the hypothesis to uniform continuity'', a claim
commonly found in the literature, is therefore a feedback-style
extrapolation (see also Section~\ref{boyer}).

The crucial point is presented by Br\aa ting in formula~(2.4) on
page~522 and the line following.  (A similar point was made by
Laugwitz \cite[p.~212]{Lau89} in 1989, in terms of the
equality~$\mu=\nu$.)  Br\aa ting documents Cauchy's use of the {\em
same\/} index~$n$, both as a subscript of a partial sum~$s_n$ of the
series~$s=\sum u_i$, and the value~$x=\frac{1}{n}$ at which the
partial sum is evaluated.  Namely, the index~$n$ appears in Br\aa
ting's formula~(2.4) as the index in an expression spelling out the
difference~$s_{n'} - s_n$, and it also appears on the next line, in
the expression~$x=\frac{1}{n}$.  In other words, Cauchy does not limit
the dynamic variable/sequential approach to his ``quantities''.%
\footnote{Cauchy's 1853 text shows that Cauchy applies such an
approach to functions, as well.  A dynamic function, such as the
sequence~$\langle s_n(x) : n\in\N \rangle$ of partial sums, is applied
by Cauchy to the quantity~$\langle \frac{1}{n}\rangle$ by evaluating
term-by-term, to obtain a new dynamic quantity~$\left\langle
s_n(\frac{1}{n}): n\in \N \right\rangle$, generating another member of
Cauchy's continuum.  While, clearly, modern constructions and concepts
such as the ultrapower construction, internal set, etc.~have no place
in Cauchy's world, a reader already familiar with the latter concepts
my find it helpful, in understanding Cauchy, to note the parallel
situation when an internal function~$[s_n]$ is applied to a
hyperreal~$[u_n]$ in a term-by-term fashion.  See also
Section~\ref{53}.}
A more detailed discussion of Cauchy's text may be found in
Section~\ref{53}.

What does emerge from Br\aa ting's analysis is that the competing
interpretations by Giusti (1984) and by Laugwitz (1987) {\em both\/}
have legitimate sources in mid-19th century work of a Swedish
mathematician who was in close contact with Cauchy, see
\cite[p.~521]{Br}.

Like L.~Carnot before him, Cauchy represented infinitesimals by null
sequences.  Cauchy spoke of variables or sequences, say~$\langle u_n:
n\in \N\rangle$, as {\em becoming\/} infinitesimals.  The precise
meaning of Cauchy's use of the verb {\em become\/} is subject to
dispute.  Meanwhile, a key question is whether, {\em after\/} becoming
an infinitesimal, such a~$\langle u_n\rangle$ is admitted to his
continuum.

We see that, in 1853, Cauchy used the expression~$x=\frac{1}{n}$ to
show that counterexamples such as Abel's 1826 ``exception'' did not
satisfy Cauchy's hypothesis.  This reveals that he {\em is\/} willing
to evaluate a function at a {\em variable\/} quantity, used as input
to the function.  The fact that Cauchy exploits such a quantity as an
{\em input\/} to his functions, suggests that quantities in the wider
sense of a B-continuum were indeed part of Cauchy's continuum, at
least at this later time.  If this interpretation is admitted, then
testing an analytical hypothesis at all members of the continuum would
naturally include testing at~$x=\frac{1}{n}$, as well.  If on the
other hand~$\langle u_n\rangle$ is not admitted as a member of the
continuum, then the continuum is restricted to what Bj\"orling called
{\em fixed\/} values (namely, belonging to an A-continuum).

\section{Cauchy's 1853 text}
\label{53}

Cauchy's text {\em Note sur les s\'eries convergentes dont les divers
termes sont des fonctions continues d'une variable r\'eelle ou
imaginaire, entre des limites donn\'ees\/} appeared in 1853, see
\cite{Ca53}.  On page~32, Cauchy recalls the definition of continuity
already mentioned in Section~\ref{two}.

Cauchy deals with a series~$s= \sum_{i=0} u_i$ with partial sum
\[
s_n=u_0 + \ldots + u_{n-1}
\]
and remainder~$r_n=s-s_n$.  He now considers~$n'>n$ and the expression
$s_{n'}-s_n= u_n+ \ldots + u_{n'-1}$, and proceeds to state his
Theorem 1 to the effect that if~$u_n$ are continuous in~$x$, and
$s_{n'}-s_n$ {\em devient toujours\/} (always becomes) infinitely
small, then the sum~$s$ will be a continuous function of the variable
$x$ \cite[p.~33]{Ca53}.

To illustrate why the series~$\sum_i \frac{\sin ix}{i}$ is not a
counterexample, he undertakes a remarkable maneuver that has sparked
controversy ever since, namely he evaluates~$s_{n'}-s_n$ at
$x=\frac{1}{n}$, with the same~$n$ appearing in the denominator of~$x$
and as a subscript in~$s_n$.  Cauchy concludes that the remainder does
not become small, by comparing it to an integral \cite[p.~34]{Ca53}.
He then proceeds to state a complex version of the same result, again
insisting on the {\em devient toujours\/} clause \cite[p.~35]{Ca53}.

It is interesting to note that in the ensuing discussion, Cauchy
evokes the property of the continuity of a function in the following
terms:
\begin{quote}
D'apr\`es ce qu'on vient de dire, une fonction monodrome de~$z$
variera par degr\'es insensibles, etc.  \cite[p.~35]{Ca53}.  
\end{quote}
The expression {\em par degr\'es insensibles\/} [by imperceptible
degrees] appears to be a reformulation of his infinitesimal definition
as stated by Cauchy on page~32.  The same expression was used by
Cauchy in his letter to Coriolis in 1837.

\section{Is the traditional reading, coherent?} 
\label{boyer}

We will build on Br\aa ting's analysis to examine a traditional
reading of Cauchy's definitions.  Cutland {\em et al.\/} note that
\begin{quote}
[Cauchy's] modification of his theorem is anything but clear if we
interpret his conception of the continuum as identical with the
`Weierstrassian' concept \cite[p.~376]{CKKR}.
\end{quote} 
We will elaborate on this comment, based on an interpretation of
Cauchy given by C.~Boyer.%
\footnote{We note in passing a curious error in \cite[p.~282]{Boy}.
Here Boyer claims that Cauchy's ``geometrical intuitions [...] led him
erroneously to believe that the continuity of a function was
sufficient [...] for the existence of a derivative.''  Boyer
continues: ``A. M. Amp\`ere also had been led by geometric
preconceptions similar to those of Cauchy to try to demonstrate the
false proposition that every continuous function has a derivative,
except for certain isolated values in the interval.''  Boyer provides
a footnote (his footnote 45) containing a reference to Jourdain
\cite{Jo}, but Jourdain's text does not bear out Boyer's claim, on the
contrary.  Jourdain makes it clear that Amp\`ere is the one who
``proved'' that every continuous function has a derivative.  Jourdain
\cite[p.~702]{Jo} discusses Amp\`ere's error in detail.  Boyer appears
to have mixed up Cauchy and Amp\`ere.  Cauchy's treatises on
differential analysis show clearly that he was aware of the fact that
possible points of non-differentiability need to be taken into account
in formulating the fundamental theorem of calculus (each point of the
first kind contributes a boundary term to the formula), and felt, as
many did in his era, that there should be only finitely many such
points.}

Boyer quotes Cauchy's definition of continuity as follows: ``the
function~$f$ is continuous within given limits if between these limits
an infinitely small increment~$i$ in the variable $x$ produces always
an infinitely small increment,~$f(x+i)-f(x)$, in the function itself''
\cite[p.~277]{Boy}.  Next, Boyer proceeds to {\em interpret\/}
Cauchy's definition of continuity as follows:
\begin{quote}
The expressions infinitely small {\em are here to be understood\/}
[...] in terms of [...]  limits: i.e., $f(x)$ is continuous within an
interval if the limit of the variable $f(x)$ as~$x$ approaches~$a$
is~$f(a)$, for any value of~$a$ within this interval
\cite[p.~277]{Boy} [emphasis added--authors]
\end{quote}

Boyer feels that infinitesimals {\em are to be understood\/} in terms
of limits.  Or perhaps they are to be understood otherwise?

Given the frequent references to Jourdain \cite{Jo} in Boyer's text,
it is worth mentioning a striking aspect of the discussion of the
notion of continuity in Jourdain \cite{Jo}: there is a total {\em
absence\/} of any claim to the effect that Cauchy based his notion of
continuity, on limits.

As we consider Boyer's interpretation in detail, we find that there
are two problems:
\begin{enumerate}
\item
historians generally agree that Cauchy did not have the notion of {\em
continuity at a point\/}.  Boyer's introduction of the {\em value\/}
$a$, and quantification over~$a$, is not present in Cauchy.  
\item
consider the function~$f(x)=\sin\frac{1}{x}$ explicitly mentioned the
1821 textbook {\em Cours d'Analyse\/} by Cauchy.  How would Cauchy
view~$f$, given his definition of continuity?
\end{enumerate}

From Boyer's post-Weierstrassian viewpoint, the function~$f$ is
continuous wherever it is defined.  However, this is not necessarily
Cauchy's view.  In Section~\ref{intro}, we followed Br\aa ting in
analyzing Cauchy's test of his condition with regard to the input~$x$
generated by the sequence
\[
\left\langle \frac{1}{n} : n\in \N \right\rangle.
\]
Now choose the infinitesimal~$i$ generated by the same sequence.  The
difference
\[
f(x+i)-f(x)=f\left(\tfrac{2}{n}\right)-f\left(\tfrac{1}{n}\right)=\sin
\tfrac{n}{2}-\sin n
\]
does not tend to zero.  If so, would~$f$ pass Cauchy's test for
continuity?

The basic problem with standard Cauchy historiography, as exemplified
by Boyer's interpretation of Cauchy's infinitesimals, seems to be as
follows.  Many historians have claimed that Cauchy modified the
hypothesis of his sum theorem in 1853, by introducing the stronger
hypothesis of uniform continuity (more precisely, uniform
convergence).  If one wishes to substantiate such a claim, then one
must interpret Cauchy's use of the term {\em always\/} as meaning that
Cauchy requires convergence not merely at the fixed numbers
(A-continuum), but also at the variable quantities (B-continuum), such
as infinitesimals.

But if one wishes to apply the Boyer infinitesimal-to-limit
translation, with an attendant interpretation of the point {\em a\/}
as a real number, then one's conceptual framework excludes the
possibility of evaluation at a variable quantity.  If one excludes
variable quantities by adhering to the infinitesimal-to-limit
translation, then one is unable to interpret Cauchy's extended
hypothesis in 1853.  If one wishes to understand Cauchy's extension of
the hypothesis, one has to jettison the automatic
translation-to-limits.  What is caught in this tightening noose is a
body of flawed Cauchy scholarship going back to Boyer or earlier.

How does the traditional approach connect Cauchy's term {\em
always\/}, to uniform convergence?  J.~L\"utzen notes that
\begin{quote}
The key word that separates [the 1853] statement from [Cauchy's]
previous [1821] statement is ``always" but only in the proof it
becomes clear what it covers \cite[p.~184]{Lut03}.
\end{quote} 
L\"utzen proceeds to reproduce a paragraph from Cauchy's proof, and
notes that Cauchy's term `` `always' covers the concept `uniform
Cauchy sequence in an interval' from which Cauchy immediately
concluded `uniform convergence in an interval'.''  L\"utzen concludes
as follows:
\begin{quote}
{\em Cauchy carefully showed\/} that a Fourier series similar to
Abel's (Cauchy did not mention Abel) does not ``always" converge in
this sense, which explains why its sum is discontinuous
\cite[p.~184]{Lut03}.  [emphasis added--authors]
\end{quote}
Now Cauchy did not use the terminology of either ``uniform Cauchy
sequence", or ``uniform convergence".  L\"utzen does not explain how
it was exactly that {\em Cauchy carefully showed\/}.  Similarly,
L\"utzen does not reproduce Cauchy's example~$x=\frac{1}{n}$ which
would have shed light on the matter, by revealing a link to a
B-continuum.

\section{Conclusion}

An examination of Cauchy's work on the sum theorem reveals that a
coherent explanation thereof requires infinitisimals to be part and
parcel of Cauchy's continuum, as they were of Leibniz's, Bernoulli's,
and Carnot's.  The historical and philosophical significance of our
analysis is the revelation that modern reception of Cauchy's
foundational work has been colored by a nominalistic attitude
resulting in an ostrich effect when it comes to appreciating Cauchy's
infinitesimals, an attitude all the more puzzling since it must
countenance an internal contradiction as analyzed in this article.
See \cite{KK11a} for a detailed examination of a nominalism inherited
from the great triumvirate.%
\footnote{\label{triumvirate}Boyer \cite[p.~298]{Boy} refers to
Cantor, Dedekind, and Weierstrass as ``the great triumvirate''.}

\appendix

\section{Spalts Kontinuum}

In a text confidently entitled ``Cauchys Kontinuum'' \cite{Sp},
D.~Spalt seeks to provide a novel interpretation of Cauchy's
foundational approach.  Spalt affirms the correctness of Cauchy's sum
theorem of 1821, and at the same time denies that Cauchy ever used
infinitesimals.

The starting point of Spalt's interpretation in \cite{Sp} is Cauchy's
double parenthesis notation.  Cauchy used such notation to signal
situations where a multiple-valuedness arises.  Spalt's interpretation
rests on the mathematical fact that if a function (a) has a closed
graph and (b) is single-valued, then it is continuous.  Was it
Cauchy's intention to define {\em continuous functions\/} in terms of
such a property?  If so, Cauchy would have called them {\em
single-valued functions\/}.  Cauchy's continuity has its source in
naive perceptual continuity.  Sensory perception experiences
continuity in terms of slight dynamic change, when an
infinitesimal~$x$-increment results in an infinitesimal change of the
dependent variable.  Having defined continuity in terms of such a rule
of transforming infinitesimals into infinitesimals both in his {\em
Cours d'Analyse\/} of 1821 and in his lectures of 1823, Cauchy again
emphasizes this point in his letter to Coriolis in 1837:~$y$ varies
{\em imperceptibly\/} with~$x$.

Whenever Cauchy used the double parenthesis notation, it is always
with reference to a {\em single\/} function~$f$, such as~$\frac{1}{x}$
or~$\sqrt{x}$ or~$\arccos x$ or~$\sin \frac{1}{x}$.  Meanwhile, Spalt
is mainly interested in applying sequences of functions to sequences
of points, as we discuss below.

The traditional interpretation as exemplified by Boyer seeks to
subsume Cauchy's infinitesimals in what Boyer sees as an inchoate
proto-Weierstrassian limit of~$f$ at, say,~$x=0$.  Spalt, meanwhile,
seeks to subsume Cauchy's infinitesimals in the calculation of the
(potentially) multiple values of~$f$ at~$x=0$.

However, Spalt's real interest is in applying the closed-graph
interpretation to Cauchy's sum theorem.  Here it is a {\em sequence\/}
of partial sums that is being evaluated at a null sequence, so as to
test the behavior of the limit at, say,~$x=0$.  There is no textual
support in Cauchy for applying the double-parenthesis notation to a
sequence.  Whenever double parentheses are used in Cauchy, it is
always with reference to a single function.

While Spalt is making a mathematically valid point that both
continuity and the sum theorem admit a ``closed graph'' interpretation
in terms of sequences, the attribution of such an interpretation to
Cauchy is not supported by textual evidence.

The following exchange, represented by individuals A and B, took place
in the fall of 2010 and illustrates well the issues involved in
evaluating Cauchy's infinitesimals.

\medskip
A. Concerning ``index": You claim that Cauchy ``suppresses the index"
in connection with his infinitely small quantities.  Where do you know
this from?  You can only suppress something you have - but Cauchy did
not have indices in connection with his infinitely small quantites.
So you insinuate these indices, but they are not Cauchyan.

\medskip
B. On page 192 of the first volume of the {\em Math.~Intelligencer\/},
Guggenheimer \cite{Gu} equips Cauchy's infinitesimals a lower index
``$n$" by writing~$\beta= \{\beta_n\}$, and acts as if Cauchy did the
same on page~26 of the {\em Cours d'Analyse\/}.  In a subsequent
issue, Gordon Fisher \cite{Fi} takes issue with this, and states that
``It is Guggenheimer who introduces the sequence~$\beta_n$ into the
definition."  Meanwhile, Cauchy does use lower indices for sequences
(though not for infinitesimals) in his proof of the intermediate value
theorem \cite[Note~III, p.~460-462]{Ca21}.

\medskip
A. Every mathematician from at least about 1750 till about 1872
thought of the points of the continuum as forming a series - but none
of them introduced this series as an indexed one.  Not Cauchy either -
whatever historians (like Guggenheimer) do say.  Cauchy's indices
(e.g. in~$\sum_i u_i(x)$) are variables with the natural numbers as
values. But as far as I know Cauchy never claimed the points of the
continuum to be capable of being indexed by the natural numbers.  This
does not prove that Cauchy doubted this possibility - but surely he
did not see any possibility of {\em how\/} to do this.

\medskip
B. Guggenheimer was certainly wrong as Fisher pointed out.  Cauchy's
variable quantities are sometimes discrete (as in 1821, when he gives
an example 1, 1/2, 1/3, etc.) and sometimes continuous (as in 1823).
As far as the continuum is concerned, Cauchy certainly never claimed
to index the points of the continuum.  But Cauchy himself gives an
example of a discrete variable quantity in 1821, while in 1823 he
worked with continuous variable quantities.  The fact that he never
labels its terms by a lower index confirms my sentiment that he wants
to underemphasize the role of this index and emphasize on the contrary
other indices, such as the order of the infinitesimal.  The use of the
term ``suppress" in this sense is legitimate whether or not the index
was there in the first place.

\medskip
A. Every renowned German mathematician of the middle of the 19th
century defined continuity initially as: infinitely small changes of
the variable produce infinitely small changes of the function.  Only
in the sequel did they give an epsilon-delta-formulation.  This
coverage includes even Weierstrass!  So there seems to be no
justification at all for later historians to claim the existence of a
conceptual difference between an A-continuum and B-continuum in the
middle of the 19th century.  Of course this judgement has to include
Cauchy as long as there is no proof of the contrary.

\medskip
B. The proof to the contrary is Cauchy's 1853 text from the middle of
the 19th century.  There are numerous other proofs, as well.  They
include Cauchy's infinitesimal definitions of ``Dirac'' delta
functions in 1827 \cite{Ca27a}.  On the other hand, what there is no
proof to at all is the idea of a ``closed graph" interpretation of
Cauchy's sum theorem \`a la Spalt.  It is true that a function is
continuous if and only if its graph is closed.  However, this has
nothing to do with the function being or not being the sum of an
infinite series.  Therefore this has nothing to do with Cauchy's sum
theorem.

\medskip
A. Could you specify where exactly such a ``proof'' is to be found?
Which are the relevant sentences?

\medskip
B. The addition of the term ``always" to the hypothesis of the sum
theorem in 1853 is interpreted by all traditional historians as adding
the condition of uniform continuity.  The meaning of the term only
becomes clear in the proof when Cauchy tests the condition at the
point $x=\frac{1}{n}$, showing that Abel's counterexample does not
satisfy the hypothesis.  One obtains uniform convergence by requiring
the remainder term to tend to zero at the points of the B-continuum in
addition to those of the A-continuum.  Cauchy tests the condition at
$x=\frac{1}{n}$.  This sequence generates an infinitesimal, i.e. a
point of a B-continuum.  L\"utzen \cite{Lut03} fails to explain this,
but Br\aa ting does.  The reference for the Dirac delta function is in
Cauchy's 1827 texts cited in Laugwitz \cite{Lau89}.  Cauchy's theory
of arbitrary {\em real\/} orders of infinitesimals is in Cauchy 1829
\cite{Ca29}, and it anticipates work on orders of growth of functions
by Stolz and du Bois-Reymond.  The latter in turn influenced Skolem's
construction of non-standard models of arithmetic.  Robinson wrote:
``It seems likely that Skolem's idea to represent infinitely large
natural numbers by number-theoretic functions which tend to infinity
(Skolem [1934]), also is related to the earlier ideas of Cauchy and du
Bois-Reymond'' \cite[p.~278]{Ro66}.

\medskip
A. Unfortunately I can't accept your ``proof".  You point to an
example and you give {\em interpretations\/}, but you don't have a
single {\em definition\/} (of Cauchy) at hand to strengthen your
position.  You claim that $1/n$ is an infinitesimal (and you mean: a
point of your B-continuum).  Cauchy himself does {\em not\/} call
$1/n$ an infinitesimal; even though he could have done, as he {\em
defines\/} an infinitesimal to be a variable converging to zero (which
$1/n$ certainly is), and so Cauchy can rely on the then common
A-continuum.

\medskip
B. The term ``always" indicates a strengthening of the hypothesis.
The hypothesis is strengthened by requiring the convergence condition
at additional members of Cauchy's continuum.  One such additional
member is generated by $1/n$.  Cauchy proceeds to use it as an input
to his functions.  This is very similar to the variable quantity
Cauchy gives as an example in 1821, namely the sequence
\[
\frac{1}{4},\; \frac{1}{3},\; \frac{1}{6},\; \frac{1}{5},\;
\frac{1}{8},\; \frac{1}{7} \;\cdots,
\]
see \cite[p.~27]{Ca21}.  These facts indicate that Cauchy {\em was\/}
working with an extended continuum.

\medskip
A.  You ground your whole thesis merely on Cauchy's term ``always"?
You are really willing to claim Cauchy to have constructed a
``B-continuum" (an outstanding mathematical construction none of
Cauchy's contemporaries ever thought of) on this single word
``always"?

\medskip
B. The B-continuum, as the name suggests, is rooted in the work of
Bernoulli.  Until the nominalistic reconstruction effected by
Weierstrass and his followers starting in 1870, most mathematicians
worked with infinitesimals and naturally envisioned an
infinitesimal-enriched continuum, where entities in addition to Stevin
(real) numbers%
\footnote{According to van der Waerden, Simon Stevin's ``general
notion of a real number was accepted, tacitly or explicitly, by all
later scientists'' \cite[p.~69]{van}.}
can be used as input to functions.  Cauchy explained the use of the
word ``always" in his {\em proof\/}.  Most historians in fact explain
the addition of the word ``always" as the addition of uniform
continuity.  This is only possible to do via a B-continuum.

\medskip
A. Your interpretation is time-dependent: Earlier than 1958 you would
and could not have given it!  But we are talking about sources from
the 19th century - and you need an interpretation which was possible
already in the 19th century, not only a century later.

\medskip
B.  You are apparently referring to the year of publication of the
work \cite{SL} by Schmieden and Laugwitz on nonstandard extensions.
But the idea that an infinitesimal is represented by a null sequence
is an ancient idea, and one that is even incorporated in the name
itself.  The word ``infinitesimal'' is a 17th century Latin formation
meaning ``infinitieth term" in a progression.  Interpreting variable
quantities as sequences is a widely accepted way of interpreting
Cauchy.  Variable quantities viewed as infinitesimals are already in
l'Hopital.  The idea that when Carnot talks about variable quantities,
he really means ``infinitesimals", seems to be widely accepted by
historians.  When Cauchy gives the same definition, should we assume
that he means something else?  The interpretation you referred to was
already possible in the 17th century.  The novelty of the 19th century
was the nominalistic transformation effected by Weierstrass that
prohibited talk about infinitesimals on pain of being declared guilty
of metaphysics.  But the ideology of the ``great triumvirate"%
\footnote{See footnote~\ref{triumvirate}.}
is being challenged by an increasingly vocal group of scholars, who
are in particular not satisfied that when Cauchy talks about
``infinitesimals", one must assume that he really means something
else.

\medskip
A.  Your B-continuum only exists (different from A-continuum) if
infinitesimals are numbers, not variables. And infinitesimal numbers
did not exist (as mathematical concepts - not as a chimera) before the
20th century.

\medskip
B. A close reading of Leibniz suggests, on the contrary, that
infinitesimals viewed as individuals/atomic entities are present in
European mathematical thinking as early as the 17th century.  As far
as the 19th century is concerned, Ehrlich~\cite{Eh} documents in
detail the development of non-Archimedean systems in Stolz, du
Bois-Reymond, and others.  The distinction between number and variable
that you insist upon is an artificial one.  The point is that Cauchy
uses infinitesimals as inputs to his functions, and operates with them
as if they were individuals/atomic entities.  In this sense they are
members of his continuum, though perhaps not of Spalts Kontinuum.


\section{Fermat, Wallis, and an ``amazingly reckless'' use of infinity}
\label{rival1}

\begin{figure}
\begin{equation*}
\xymatrix@C=95pt{{} \ar@{-}[rr] \ar@{-}@<-0.5pt>[rr]
\ar@{-}@<0.5pt>[rr] & {} \ar@{->}[d]^{\hbox{st}} & \hbox{\quad
B-continuum} \\ {} \ar@{-}[rr] & {} & \hbox{\quad A-continuum} }
\end{equation*}
\caption{\textsf{Thick-to-thin: taking standard part (the thickness of
the top line is merely conventional, and meant to suggest the presence
of additional numbers, such as infinitesimals)}}
\label{31}
\end{figure}

A Leibnizian definition of the derivative as the infinitesimal
quotient
\begin{equation*}
\frac{\Delta y}{\Delta x},
\end{equation*} 
whose logical weakness was criticized by Berkeley, was modified by
A.~Robinson by exploiting a map called {\em the standard part\/},
denoted~``st'', from the finite part of a B-continuum (for
``Bernoullian''), to the A-continuum (for ``Archimedean''), as
illustrated in Figure~\ref{31}.%
\footnote{In the context of the hyperreal extension of the real
numbers, the map ``st'' sends each finite point~$x$ to the real point
st$(x)\in \R$ infinitely close to~$x$.  In other words, the map ``st''
collapses the cluster (halo) of points infinitely close to a real
number~$x$, back to~$x$.}
Here two points of a B-continuum have the same image under ``st'' if
and only if they are equal up to an infinitesimal.  

This section analyzes the historical seeds of Robinson's theory, in
the work of Fermat, Wallis, as well as Barrow.%
\footnote{\label{barrow}While Barrow's role is also critical, we will
mostly concentrate on Fermat and Wallis.}
The key concept here is that of {\em adequality\/} (see below).  It
should be kept in mind that Fermat never considered the local slope of
a curve.  Therefore one has to be careful not to attribute to Fermat
mathematical content that could not be there.  On the other hand,
Barrow did study curves and their slope.  Furthermore, Barrow
exploited Fermat's adequality in his work~\cite[p.~252]{Barr}, as
documented by H.~Breger \cite[p.~198]{Bre94}.

The binary relation of ``equality up to an infinitesimal'' was
anticipated in the work of Pierre de Fermat.  Fermat used a term
usually translated into English as ``adequality''.%
\footnote{In French one uses {\em ad\'egalit\'e, ad\'egal\/}, see
\cite[p.~73]{Je}.}
Andr\'e Weil writes as follows:
\begin{quote} 
Fermat [...] developed a method which slowly but surely brought him
very close to modern infinitesimal concepts.  What he did was to write
congruences between functions of~$x$ modulo suitable powers
of~$x-x_0$; for such congruences, he introduces the technical term
{\em adaequalitas, adaequare\/}, etc., which he says he has borrowed
from Diophantus.  As Diophantus V.11 shows, it means an approximate
equality, and this is indeed how Fermat explains the word in one of
his later writings \cite[p.~1146]{We}.
\end{quote} 
Weil \cite[p.~1146, footnote~5]{We} then supplies the following quote
from Fermat:
\begin{quote}
{\em Adaequetur, ut ait Diophantus,%
\footnote{The original term in Diophantus is~$\pi \alpha \rho \iota
\sigma \mbox{\it\!\'o} \tau \eta \varsigma$, see Weil
\cite[p.~28]{We84}.}
aut fere aequetur\/}; in Mr.~Mahoney's translation: ``adequal, or
almost equal" (p. 246).
\end{quote}
Here Weil is citing Mahoney \cite[p.~246]{Mah73}
(cf. \cite[p.~247]{Mah94}).  Mahoney similarly mentions the meaning of
``approximate equality'' or ``equality in the limiting case'' in
\cite[p.~164, end of footnote~46]{Mah73}.  Mahoney also points out
that the term ``adequality'' in Fermat has additional meanings.  The
latter are emphasized in a recent text by E.~Giusti~\cite{Giu}, who is
sharply critical of Breger~\cite{Bre94}.  While the review~\cite{We}
by Weil is similarly sharply critical of Mahoney, both agree that the
meaning of ``approximate equality'', leading into infinitesimal
calculus, is at least {\em one of the meanings\/} of the term {\em
adequality\/} for Fermat.%
\footnote{Jensen similarly describes adequality as approximate
equality, and describes neglected terms as {\em infinitesimals\/} in
\cite[p.~82]{Je}.  Struik notes that ``Fermat uses the term to denote
what we call a limiting process'' \cite[p.~220, footnote~5]{Stru}.
K.~Barner~\cite{Ba} compiled a useful bibliography on Fermat's
adequality, including many authors we have not mentioned here.}

\begin{figure}
\includegraphics[height=2in]{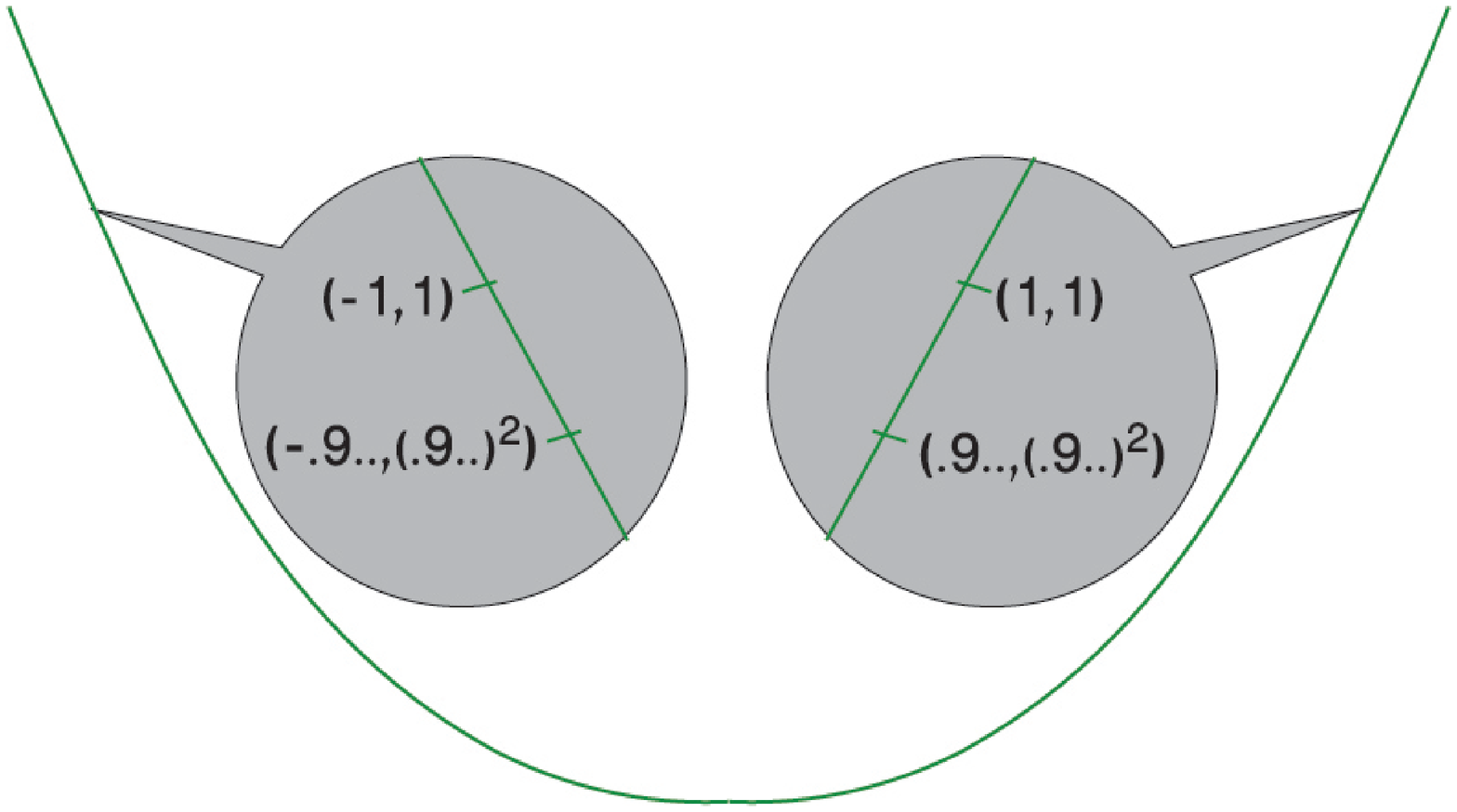}
\caption{\textsf{Differentiating~$y=f(x)=x^2$ at~$x=1$ yields
$\tfrac{\Delta y}{\Delta x} = \tfrac{f(.9..) - f(1)}{.9..-1} =
\tfrac{(.9..)^2 - 1}{.9..-1} = \tfrac{(.9.. - 1)(.9.. + 1)}{.9..-1} =
.9.. + 1 \approx 2.$ Here~$\approx$ is the relation of being
infinitely close (adequal).  Hyperreals of the form~$.9..$ are
discussed in \cite{KK2}}}
\label{jul10}
\end{figure}

\begin{figure}
\includegraphics[height=2.3in]{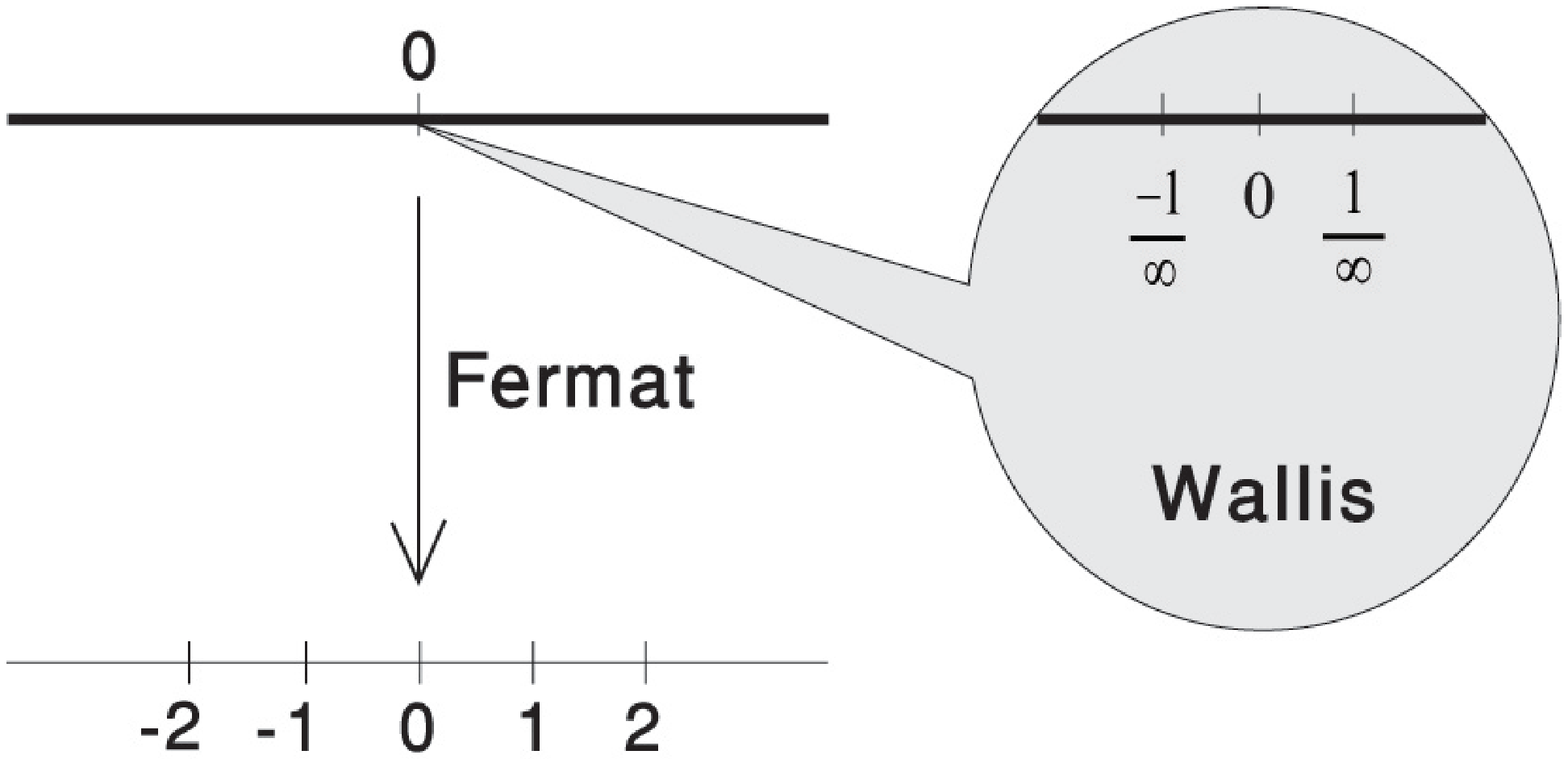}
\caption{\textsf{Zooming in on Wallis's
infinitesimal~$\frac{1}{\infty}$, which is adequal to~$0$ in Fermat's
terminology}}
\label{FermatWallis}
\end{figure}

This meaning was aptly summarized by J.~Stillwell.  Stillwell's
historical presentation is somewhat simplified, and does not
sufficiently distinguish between the seeds actually present in Fermat,
on the one hand, and a modern interpretation thereof, on the other,%
\footnote{See main text around footnote~\ref{barrow} above for a
mention of Barrow's role, documented by H~Breger.}
but he does a splendid job of explaining the mathematical background
for the uninitiated.  Thus, he notes that~$2x+dx$ is not equal to~$2x$
(see Figure~\ref{jul10}), and writes:
\begin{quote}
Instead, the two are connected by a looser notion than equality that
Fermat called adequality.  If we denote adequality by~$=_{ad}$, then
it is accurate to say that
\[
2x+dx=_{ad}2x,
\]
and hence that~$dy/dx$ for the parabola is adequal to~$2x$.
Meanwhile,~$2x+dx$ is not a number, so~$2x$ is the only number to
which~$dy/dx$ is adequal.  This is the true sense in which~$dy/dx$
represents the slope of the curve~\cite[p.~91]{Sti}.
\end{quote}
Stillwell points out that
\begin{quote}
Fermat introduced the idea of adequality in 1630s but he was ahead of
his time.  His successors were unwilling to give up the convenience of
ordinary equations, preferring to use equality loosely rather than to
use adequality accurately.  The idea of adequality was revived only in
the twentieth century, in the so-called non-standard analysis
\cite[p.~91]{Sti}.
\end{quote}
We will refer to the map from the (finite part of the) B-continuum to
the A-continuum as the Fermat-Robinson standard part, see
Figure~\ref{FermatWallis}.

As far as the logical criticism formulated by Rev.~George is
concerned, Fermat's adequality had pre-emptively provided the seeds of
an answer, a century before the bishop ever lifted up his pen to write
{\em The Analyst\/}~\cite{Be}.

\begin{figure}
\includegraphics[height=2.3in]{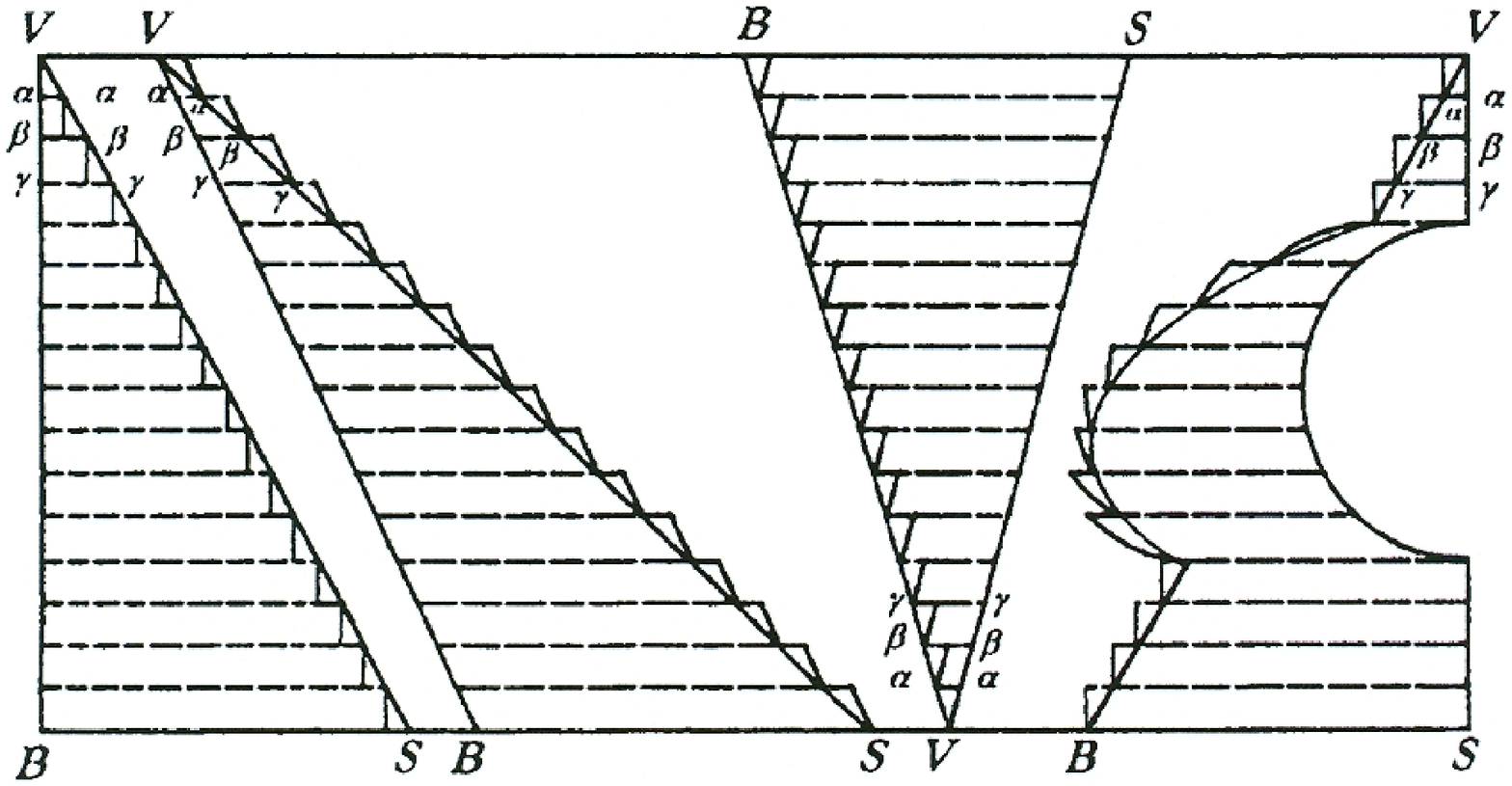}
\caption{\textsf{Area calculations in Wallis: slicing it up into
dilatable parallelograms of infinitesimal altitude}}
\label{Wallis}
\end{figure}

Fermat's contemporary John Wallis, in a departure from Cavalieri's
focus on the geometry of indivisibles, emphasized the arithmetic of
infinitesimals, see J.~Stedall's introduction in \cite{Wa}.  To
Cavalieri, a plane figure is made up of lines; to Wallis, it is made
of parallelograms of infinitesimal altitude.  Wallis transforms this
insight into symbolic algebra over the~$\infty$ symbol which he
introduced.  He exploits formulas
like~$\infty\times\frac{1}{\infty}=1$ in his calculations of areas.
Thus, in proposition 182 of his {\em Arithmetica Infinitorum\/},
Wallis partitions a triangle of altitude~$A$ and base~$B$ into a
precise number~$\infty$ of ``parallelograms" of infinitesimal
width~$\frac{A}{\infty}$, see Figure~\ref{Wallis} (copied from
\cite[p.~170]{Mun}).

He then computes the combined length of the bases of the
parallelograms to be~$\frac{B}{2}\infty$, and finds the area to be
\begin{equation}
\label{wallis}
\frac{A}{\infty}\times \frac{B}{2}\infty = \frac{AB}{2}.
\end{equation}
Wallis used an actual infinitesimal~$\frac{1}{\infty}$ in
calculations as if it were an ordinary number, anticipating Leibniz's
law of continuity.

Wallis's area calculation~\eqref{wallis} is reproduced by J.~Scott,
who notes that Wallis
\begin{quote}
treats infinity as though the ordinary rules of arithmetic could be
applied to it \cite[p.~20]{Sco}.
\end{quote}
Such a treatment of infinity strikes Scott as something of a blemish,
as he writes:
\begin{quote}
But this is perhaps understandable.  For many years to come the
greatest confusion regarding these terms persisted, and even in the
next century they continued to be used in what appears to us an
amazingly reckless fashion \cite[p.~21]{Sco}.
\end{quote}

What is the source of Scott's confidence in dismissing Wallis's use of
infinity as ``reckless''?  Scott identifies it on the preceding page
of his book; it is, predictably, the triumvirate ``modern conception
of infinity'' \cite[p.~19]{Sco}.  Scott's tunnel A-continuum vision
blinds him to the potential of Wallis's vision of infinity.  But this
is perhaps understandable.  Many years separate Scott from Robinson's
theory which in particular empowers Wallis's calculation.  The lesson
of Scott's condescending steamrolling of Wallis's infinitesimal
calculation could be taken to heart by historians who until this day
cling to a nominalistic belief that Robinson's theory has little
relevance to the history of mathematics in the 17th century.

\section{Rival continua}
\label{rival2}

This section summarizes a 20th century implementation of the
B-continuum, not to be confused with incipient notions of such a
continuum found in earlier centuries.  An alternative implementation
has been pursued by Lawvere, John L. Bell \cite{Bel08, Bel}, and
others.

We illustrate the construction by means of an infinite-resolution
microscope in Figure~\ref{FermatWallis}.  We will denote such a
B-continuum by the new symbol \RRR{} (``thick-R'').  Such a continuum
is constructed in formula~\eqref{bee}.  We will also denote its finite
part, by
\begin{equation*}
\RRR_{<\infty} = \left\{ x\in \RRR : \; |x|<\infty \right\},
\end{equation*}
so that we have a disjoint union
\begin{equation}
\RRR= \RRR_{<\infty} \cup \RRR_{\infty},
\end{equation}
where~$\RRR_{\infty}$ consists of unlimited hyperreals (i.e., inverses
of nonzero infinitesimals).

The map ``st'' sends each finite point~$x\in \RRR$, to the real point
st$(x)\in \R$ infinitely close to~$x$, as follows:%
\footnote{This is the Fermat-Robinson standard part whose seeds in
Fermat's adequality were discussed in Appendix~\ref{rival1}.}
\begin{equation*}
\xymatrix{\quad \RRR_{{<\infty}}^{~} \ar[d]^{{\rm st}} \\ \R}
\end{equation*}
Robinson's answer to Berkeley's {\em logical criticism\/} (see
D.~Sherry \cite{She87}) is to define the derivative as
\begin{equation*}
\hbox{st} \left( \frac{\Delta y}{\Delta x} \right),
\end{equation*}
instead of~$\Delta y/\Delta x$.

Note that both the term ``hyper-real field'', and an ultrapower
construction thereof, are due to E.~Hewitt in 1948, see
\cite[p.~74]{Hew}.  In 1966, Robinson referred to the 
\begin{quote}
theory of hyperreal fields (Hewitt [1948]) which ... can serve as
non-standard models of analysis \cite[p.~278]{Ro66}.
\end{quote}
The {\em transfer principle\/} is a precise implementation of
Leibniz's heuristic {\em law of continuity\/}: ``what succeeds for the
finite numbers succeeds also for the infinite numbers and vice
versa'', see~\cite[p.~266]{Ro66}.  The transfer principle, allowing an
extention of every first-order real statement to the hyperreals, is a
consequence of the theorem of J.~{\L}o{\'s} in 1955, see~\cite{Lo},
and can therefore be referred to as a Leibniz-{\L}o{\'s} transfer
principle.  A Hewitt-{\L}o{\'s} framework allows one to work in a
B-continuum satisfying the transfer principle.  To elaborate on the
ultrapower construction of the hyperreals, let~$\Q^\N$ denote the ring
of sequences of rational numbers.  Let
\begin{equation*}
\left( \Q^\N \right)_C
\end{equation*}
denote the subspace consisting of Cauchy sequences.  The reals are by
definition the quotient field
\begin{equation}
\label{real}
\R:= \left. \left( \Q^\N \right)_C \right/ \mathcal{F}_{\!n\!u\!l\!l},
\end{equation}
where~$\mathcal{F}_{\!n\!u\!l\!l}$ contains all null sequences.
Meanwhile, an infinitesimal-enriched field extension of~$\Q$ may be
obtained by forming the quotient
\begin{equation*}
\left.  \Q^\N \right/ \mathcal{F}_{u}.
\end{equation*}
Here a sequence~$\langle u_n : n\in \N \rangle$ is
in~$\mathcal{F}_{u}$ if and only if the set of indices
\[
\{ n \in \N : u_n = 0 \}
\]
is a member of a fixed ultrafilter.%
\footnote{In this construction, every null sequence defines an
infinitesimal, but the converse is not necessarily true.  Modulo
suitable foundational material, one can ensure that every
infinitesimal is represented by a null sequence; an appropriate
ultrafilter (called a {\em P-point\/}) will exist if one assumes the
continuum hypothesis, or even the weaker Martin's axiom.  See Cutland
{\em et al\/} \cite{CKKR} for details.}
See Figure~\ref{helpful}.

\begin{figure}
\begin{equation*}
\xymatrix{ && \left( \left. \Q^\N \right/ \mathcal{F}_{\!u}
\right)_{<\infty} \ar@{^{(}->} [rr]^{} \ar@{->>}[d]^{\rm st} &&
\RRR_{<\infty} \ar@{->>}[d]^{\rm st} \\ \Q \ar[rr] \ar@{^{(}->} [urr]
&& \R \ar[rr]^{\simeq} && \R }
\end{equation*}
\caption{\textsf{An intermediate field~$\left. \Q^\N \right/
\mathcal{F}_{\!u}$ is built directly out of~$\Q$}}
\label{helpful}
\end{figure}

To give an example, the sequence
\begin{equation}
\label{infinitesimal}
\left\langle \tfrac{(-1)^n}{n} \right\rangle
\end{equation}
represents a nonzero infinitesimal, whose sign depends on whether or
not the set~$2\N$ is a member of the ultrafilter.  To obtain a full
hyperreal field, we replace~$\Q$ by~$\R$ in the construction, and form
a similar quotient
\begin{equation}
\label{bee}
\RRR:= \left.  \R^\N \right/ \mathcal{F}_{u}.
\end{equation}
We wish to emphasize the analogy with formula~\eqref{real} defining
the A-continuum.  Note that, while the leftmost vertical arrow in
Figure~\ref{helpful} is surjective, we have
\begin{equation*}
\left( \Q^\N / \mathcal{F}_{u} \right) \cap \R = \Q.
\end{equation*}
A more detailed discussion of this construction can be found in the
book by M.~Davis~\cite{Dav77}.  
%
%
See also B\l aszczyk \cite{Bl} for some philosophical implications.
More advanced properties of the hyperreals such as saturation were
proved later, see Keisler \cite{Kei} for a historical outline.  A
helpful ``semicolon'' notation for presenting an extended decimal
expansion of a hyperreal was described by A.~H.~Lightstone~\cite{Li}.
See also P.~Roquette \cite{Roq} for infinitesimal reminiscences.  A
discussion of infinitesimal optics is in K.~Stroyan \cite{Str},
H.~J.~Keisler~\cite{Ke}, D.~Tall~\cite{Ta80}, L.~Magnani \&
R.~Dossena~\cite{MD, DM}, and Bair \& Henry \cite{BH}.

Applications of the B-continuum range from aid in teaching calculus
\cite{El, KK1, KK2, Ta91, Ta09a} to the Bolzmann equation (see
L.~Arkeryd~\cite{Ar81, Ar05}); modeling of timed systems in computer
science (see H.~Rust \cite{Rust}); mathematical economics (see
R.~Anderson \cite{An00}); mathematical physics (see Albeverio {\em et
al.\/} \cite{Alb}); etc.

\section*{Acknowledgments}

The authors are grateful to the anonymous referee for perceptive
comments that helped improve an earlier version of the manuscript.
Hilton Kramer's influence is obvious throughout.

\medskip\noindent
{\bf Karin Usadi Katz} has taught mathematics at Michlelet Banot
Lustig, Ramat Gan, Israel.  Two of her joint studies with Mikhail Katz
were published in {\em Foundations of Science\/}: ``A Burgessian
critique of nominalistic tendencies in contemporary mathematics and
its historiography" and ``Stevin numbers and reality", online
respectively at

http://dx.doi.org/10.1007/s10699-011-9223-1 and at

http://dx.doi.org/10.1007/s10699-011-9228-9

A joint study with Mikhail Katz entitled ``Meaning in classical
mathematics: is it at odds with Intuitionism?" is due to appear in
{\em Intellectica\/}.

\medskip\noindent
{\bf Mikhail G. Katz} is Professor of Mathematics at Bar Ilan
University, Ramat Gan, Israel.  Two of his joint studies with Karin
Katz were published in {\em Foundations of Science\/}: ``A Burgessian
critique of nominalistic tendencies in contemporary mathematics and
its historiography" and ``Stevin numbers and reality", online
respectively at

http://dx.doi.org/10.1007/s10699-011-9223-1 and at

http://dx.doi.org/10.1007/s10699-011-9228-9

A joint study with Karin Katz entitled ``Meaning in classical
mathematics: is it at odds with Intuitionism?" is due to appear in
{\em Intellectica\/}.

A joint study with Alexandre Borovik entitled ``Who gave you the
Cauchy--Weierstrass tale?  The dual history of rigorous calculus''
appeared in {\em Foundations of Science\/}, online at

http://dx.doi.org/10.1007/s10699-011-9235-x

A joint study with David Tall, entitled ``The tension between
intuitive infinitesimals and formal mathematical analysis", is due to
appear as a chapter in a book edited by Bharath Sriraman, see

\noindent
http://www.infoagepub.com/products/Crossroads-in-the-History-of-Mathematics

\end{document}